\documentclass[11pt]{article}
\usepackage{amssymb,latexsym,amsmath,amsbsy,amsthm,cite}
\newtheorem{theo}{Theorem}
\def\nn{\nonumber}
\def\ssl{\mathfrak{sl}}
\def\g{\mathfrak{g}}

\begin{document}
\begin{center}
{\Large \bf
The exceptional Lie algebra $\g_2$ generated by three generators\\[2mm]
subject to quadruple relations} \\[5mm]
{\bf N.I.~Stoilova}\footnote{E-mail: stoilova@inrne.bas.bg}\\[1mm] 
Institute for Nuclear Research and Nuclear Energy, Bulgarian Academy of Sciences,\\ 
Boul.\ Tsarigradsko Chaussee 72, 1784 Sofia, Bulgaria\\[2mm] 
{\bf J.\ Van der Jeugt}\footnote{E-mail: Joris.VanderJeugt@UGent.be}\\[1mm]
Department of Applied Mathematics, Computer Science and Statistics, Ghent University,\\
Krijgslaan 281-S9, B-9000 Gent, Belgium.
\end{center}

\vskip 2 cm

\begin{abstract}
\noindent 
In this short communication we show how the Lie algebra $\g_2$ can easily be described as a free Lie algebra on 3 generators, 
subject to some simple quadruple relations for these generators.
\end{abstract}

\vskip 10mm

\setcounter{equation}{0}

The exceptional Lie groups and algebras and their descriptions have always attracted much attention.
Many presentations and realizations can be found in the literature and in books~\cite{Jacobson,Jacobson2,Humphreys,Fulton}, especially for the Lie algebra $\g_2$ (often also denoted by $G_2$).
The Lie group $G_2$ is commonly defined as the group of automorphisms of the octonions~\cite{Baez}, 
and thus the Lie algebra $\g_2$ is identified as the derivation algebra of octonions.
There are also realizations of $\g_2$ as a derivation algebra of other non-associative algebras~\cite{Basak, Elduque},
combinatorial constructions~\cite{Wildberger} or constructions starting from different automorphism groups~\cite{Wilson}.
For a survey, see~\cite{Draper} or~\cite{Agricola}.

The common description of $\g_2$ as a Lie algebra with generators and relations is in terms of the Chevalley generators and the Serre relations~\cite{Humphreys}.
As far as we know a presentation of $\g_2$ as a Lie algebra with 3 generators subject to some quadruple relations has not been given; 
at least, we were not able to find it in the literature.
Therefore we think it is worthwhile to correspond our simple result to the mathematics community.

\begin{theo}
The Lie algebra $\g$ on three generators $x_1, x_2, x_3$ subject to the following list of quadruple relations is equal to $\g_2$:
\begin{enumerate}
\item $[x_i,[x_j,[x_i,x_k]]] = 2\epsilon_{ijk} x_i$,
\item $[x_i,[x_i,[x_j,x_k]]] = 4\epsilon_{ijk} x_i$,
\item $[x_i,[x_j,[x_j,x_k]]] = 6\epsilon_{ijk} x_j$.
\end{enumerate}
Herein, $\epsilon_{ijk}$ is the common Levi-Civita symbol in three dimensions: 
$\epsilon_{ijk}$ is 1 if $(i,j,k)$ is an even permutation of $(1,2,3)$, $-1$ if it is an odd permutation of $(1,2,3)$, and 0 otherwise.
\end{theo}

Note that this list describes all possible quadruple relations for four elements with three indices from $\{1,2,3\}$, since among 4 indices at least two should be equal.

The most straightforward proof of this result is by identifying the independent elements of $\g$, and computing their commutation relations.

Let us first count the number of independent elements of $\g$. 
There are three elements of degree~1 in the generators: $x_1$, $x_2$ and $x_3$.
Using the anti-symmetry of the Lie bracket, there are also three elements of degree~2 in the generators: $[x_2,x_3]$, $[x_1,x_3]$ and $[x_1,x_2]$.
Next, the elements of degree~3 are of the form $[x_i,[x_j,x_k]]$, where $i\in\{1,2,3\}$ and $(j,k)\in\{ (2,3),(1,3),(1,2)\}$. 
However, among these 9 elements there is one linear relation following from the Jacobi identity, namely
\[
[x_1,[x_2,x_3]]+[x_2,[x_3,x_1]]+[x_3,[x_1,x_2]]=0.
\]
This implies that there are 8 linearly independent elements of degree~3 in the generators.
The quadruple relations imply that there are no independent elements of degree~4 or higher.
So all together, $\g$ has dimension~14.

In order to identify $\g$ with $\g_2$ it is sufficient to show that $\g$ is a simple Lie algebra.
We shall do more, and give the complete table of brackets among the 14 basis elements (from which it also follows that $[\g,\g]=\g$).
For this purpose, let us introduce some further notation.
The elements of degree~2 are denoted by:
\begin{equation}
y_1=\frac12 [x_2,x_3],\qquad y_2=\frac12 [x_3,x_1],\qquad y_3=\frac12 [x_1,x_2].
\label{y}
\end{equation}
The elements of degree~3 in the generators are denoted by:
\begin{align}
& a_{12}=\frac13 [x_2,y_1], \qquad a_{23}=\frac13 [x_3,y_2],\qquad a_{13}=\frac13 [x_3,y_1], \nn \\ 
& a_{21}=\frac13 [x_1,y_2], \qquad a_{32}=\frac13 [x_2,y_3],\qquad a_{31}=\frac13 [x_1,y_3], \label{a} \\
& h_1=\frac13([x_1,y_1]-[x_2,y_2]), \qquad h_2=\frac13([x_2,y_2]-[x_3,y_3]). \nn
\end{align}
It is now a simple matter to determine the table of brackets between these elements, using the definitions~\eqref{y}-\eqref{a}, the quadruple relations 1--3, and the Jacobi identity to rewrite elements of degree~4 in the form of relations 1--3.

For example, let us compute a bracket corresponding to an element of degree 5 in the generators:
\begin{align}
[y_1,a_{23}]&= \frac12 [[x_2,x_3],a_{23}] = \frac12 [x_2,[x_3,a_{23}]] - \frac12 [x_3,[x_2,a_{23}]] \nn \\
&= \frac1{12} [x_2,[x_3, [x_3,[x_3,x_1]]]] - \frac1{12} [x_3,[x_2,[x_3,[x_3,x_1]]]] = \frac1{12}[x_2,0] -\frac1{12} [x_3, 6 x_3] =0. \label{d5}
\end{align}
And as a second example, we compute an element of degree 6:
\begin{align}
[a_{12},a_{23}]&= \frac13 [[x_2,y_1],a_{23}]= \frac13 [x_2,[y_1,a_{23}]]-\frac13 [y_1,[x_2,a_{23}]] \nn\\
&= 0 -\frac13 [y_1, \frac16 [x_2,[x_3,[x_3,x_1]]] ] = -\frac13 [y_1,x_3]=a_{13}, \label{d6}
\end{align}
where the first term vanishes using the previous computation~\eqref{d5}.
Using such manipulations, the complete commutator table (Table~1) is computed.

\begin{table}[tbh]
\addtolength{\tabcolsep}{-2pt}  
\begin{center}
\begin{tabular}{c||c|c|c|c|c| c|c|c| c|c|c| c|c|c|}
 $[\cdot,\cdot]$    & $h_1$ & $h_2$ & $a_{12}$ & $a_{13}$ & $a_{23}$ & $a_{21}$ & $a_{31}$ & $a_{32}$ & $x_1$ & $x_2$ & $x_3$ & $y_1$ & $y_2$ & $y_3$ \\ \hline\hline
$h_1$ & $0$ & $0$ & $2a_{12}$ & $a_{13}$ & $-a_{23}$ & $-2a_{21}$ & $-a_{31}$ & $a_{32}$ & $-x_1$ & $x_2$ & 0 & $y_1$ & $-y_2$ & 0 \\ \hline
$h_2$ &  & $0$ & $-a_{12}$ & $a_{13}$ & $2a_{23}$ & $a_{21}$ & $-a_{31}$ & $-2a_{32}$ & 0 & $-x_2$ & $x_3$ & 0 & $y_2$ & $-y_3$ \\ \hline
$a_{12}$ &  &  &  0 & 0 & $a_{13}$ & $h_1$ & $-a_{32}$ & 0 & $-x_2$ & 0 & 0 & 0 & $y_1$ & 0 \\ \hline
$a_{13}$ &  &  &  & 0 & 0 & $-a_{23}$ & $h_1\!+\!h_2$ & $a_{12}$ & $-x_3$ & 0 & 0 & 0 & 0 & $y_1$ \\ \hline
$a_{23}$ &  &  &  &  & 0 & 0 & $a_{21}$ & $h_2$ & 0 & $-x_3$ & 0 & 0 & 0 & $y_2$ \\ \hline
$a_{21}$ &  &  &  &  &  & 0 & 0 & $-a_{31}$ & 0 & $-x_1$ & 0 & $y_2$ & 0 & 0 \\ \hline
$a_{31}$ &  &  &  &  &  &  & 0 & 0 & 0 & 0 & $-x_1$ & $y_3$ & 0 & 0 \\ \hline
$a_{32}$ &  &  &  &  &  &  &  & 0 & 0 & 0 &$-x_2$ & 0 & $y_3$ & 0 \\ \hline
$x_1$    &  &  &  &  &  &  &  &  & 0 & $2y_3$ & $-2y_2$ & $2h_1\!+\!h_2$ & $3a_{21}$ & $3a_{31}$ \\ \hline
$x_2$    &  &  &  &  &  &  &  &  &  & 0 & $2y_1$ & $3a_{12}$ & $-h_1\!+\!h_2$ & $3a_{32}$ \\ \hline
$x_3$    &  &  &  &  &  &  &  &  &  &  & 0 & $3a_{13}$ & $3a_{23}$ & $-h_1\!-\!2h_2$ \\ \hline
$y_1$    &  &  &  &  &  &  &  &  &  &  &  & 0 & $2x_3$ & $-2x_2$ \\ \hline
$y_2$    &  &  &  &  &  &  &  &  &  &  &  &  & 0 & $2x_1$ \\ \hline
$y_3$    &  &  &  &  &  &  &  &  &  &  &  &  &  & 0 \\ \hline
\end{tabular}
\end{center}
\addtolength{\tabcolsep}{2pt}  
\caption{Table of brackets among the 14 basis elements}
\end{table}

From the commutator table, the subalgebra structure of the basis is obvious. 
Clearly, the elements $h_i$ and $a_{ij}$ satisfy the standard commutation relations of $\ssl(3)$:
in the defining representation of $\ssl(3)$ in terms of $3\times 3$-matrices, the elements $a_{ij}$ 
can be realized as $e_{ij}$ (a matrix with 1 on position $(i,j)$ and zeros elsewhere) and $h_i$ as $e_{ii}-e_{i+1,i+1}$.
The elements $x_1,x_2,x_3$ are an $\ssl(3)$ triple, and $y_1,y_2,y_3$ a dual $\ssl(3)$ triple.

\section*{Acknowledgments}
Both authors were supported by the Bulgarian National Science Fund, grant KP-06-N28/6.

\end{document}